\documentclass{article}

\usepackage{arxiv}

\usepackage[utf8]{inputenc} % allow utf-8 input
\usepackage[T1]{fontenc}    % use 8-bit T1 fonts
\usepackage{hyperref}       % hyperlinks
\usepackage{url}            % simple URL typesetting
\usepackage{booktabs}       % professional-quality tables
\usepackage{amsfonts}       % blackboard math symbols
\usepackage{nicefrac}       % compact symbols for 1/2, etc.
\usepackage{microtype}      % microtypography
\usepackage{lipsum}		% Can be removed after putting your text content
\usepackage{graphicx}

\usepackage{moreverb}

\usepackage{graphicx, subfigure}

\newcommand\BibTeX{{\rmfamily B\kern-.05em \textsc{i\kern-.025em b}\kern-.08em
T\kern-.1667em\lower.7ex\hbox{E}\kern-.125emX}}

\usepackage{graphicx}
\RequirePackage{fix-cm}
\usepackage{hyperref}
 \usepackage{epsfig} % for postscript graphics files
 \usepackage{amssymb}
 \usepackage{amsmath}
\usepackage{multirow}
\usepackage{epsfig} % for postscript graphics files
\input{epsf}
\usepackage{graphicx}
\usepackage{graphicx, subfigure}
\usepackage[utf8]{inputenc}
\usepackage[english]{babel}
\usepackage{amsmath}
\usepackage{graphicx}      % include this line if your document contains figures
\usepackage{natbib}
\usepackage{graphics} % for pdf, bitmapped graphics files
\usepackage{epsfig} % for postscript graphics files
\usepackage{mathptmx} % assumes new font selection scheme installed
\usepackage{times} % assumes new font selection scheme installed
\usepackage{amsmath} % assumes amsmath package installed
\usepackage{amssymb}  % assumes amsmath package installed
\usepackage{multirow}
\input{epsf}
\usepackage{graphicx}
\usepackage{graphicx, subfigure}
\usepackage[utf8]{inputenc}
\usepackage[english]{babel}
\usepackage{matlab-prettifier}
\usepackage{listings}% http://ctan.org/pkg/listings
\usepackage[inline]{enumitem}

\def\beq{\begin{equation}}
\def\eeq{\end{equation}}
\def\baq{\begin{eqnarray}}
\def\eaq{\end{eqnarray}}
\def\bal{\begin{align} }
\def\eal{\end{align} }
\def\bc{\begin{center}}
\def\ec{\end{center}}

\def\gr{\gamma_{r}}
\def\gr1{\gamma_{r1}}

\title{Automatic Generation of Examinations in the Automatic Control Courses }
\author{ Alexander Stotsky \&  Torsten Wik \\
         Systems \& Control \\
    Department of Electrical Engineering \\
     Chalmers University of Technology  \\
     Gothenburg  SE - 412 96, Sweden  \\
	\texttt{alexander.stotsky@chalmers.se} \\
        \texttt{torsten.wik@chalmers.se  } }
\hypersetup{
pdftitle={A template for the arxiv style},
pdfsubject={q-bio.NC, q-bio.QM},
pdfauthor={David S.~Hippocampus, Elias D.~Striatum},
pdfkeywords={First keyword, Second keyword, More},
}
\date{}
\begin{document}
\maketitle

\begin{abstract}
~~Final written examination is the most important part of summative
assessment in automatic control courses. Preparation of the examinations with a given number
of points according to the concept of Constructive Alignment (which could be the main concept in future automatic control education) takes significant amount of time of the educator and motivates development of a toolkit for automatic compilation of examination problems.
A decision support Matlab/\LaTeX\ toolkit based on random number generators for selection of examination problems
is described in this report to facilitate the alignment.
The toolkit allows application of Stepwise Constructive Alignment (a new method described in this report), where the alignment is achieved by a number of software runs associated with random trials. In each step the educator manually selects suitable problems before each run based on evaluation of the random choice from the previous run.
Automatic generation of the examination together with solutions for the course  'Process control and measurement techniques' is presented as an example.
\end{abstract}


\maketitle

%===============================================================================
\begin{figure*}[!ht]
\begin{center}
  \includegraphics[width=16cm]{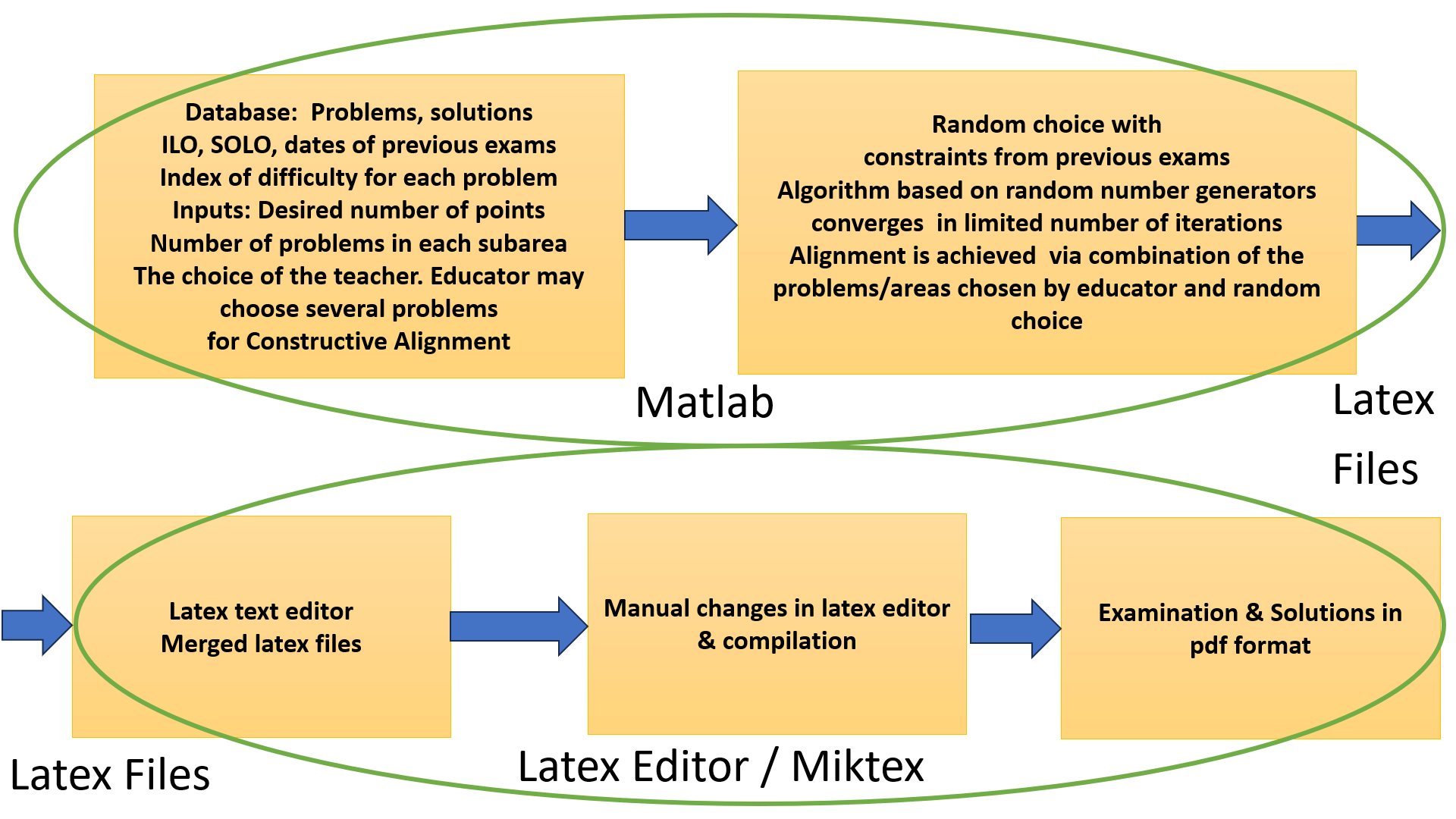}
\end{center}
\caption{\small{The flow chart for automatic generation of examination.
 This system consists of two
large parts. The first part is written in Matlab and is associated with the database, which contains examination problems and solutions, learning outcomes, SOLO model, indices of difficulty and  dates of previous examinations. The educator specifies
the number of points and the number of problems in each area to be included in the examination. The educator may select several problems manually for CA.
The Matlab program creates the \LaTeX\ file by merging items of randomly selected problems in the subareas.
The merged \LaTeX\ file can be reviewed in WinEdt editor for example, where the educator gets opportunity to edit the examination and solutions.
The \LaTeX\ source file is compiled by Miktex which produces the pdf output with examination and solutions.
    }}

\label{figs}
\end{figure*}
\begin{figure}
\begin{center}
  \includegraphics[width=9cm]{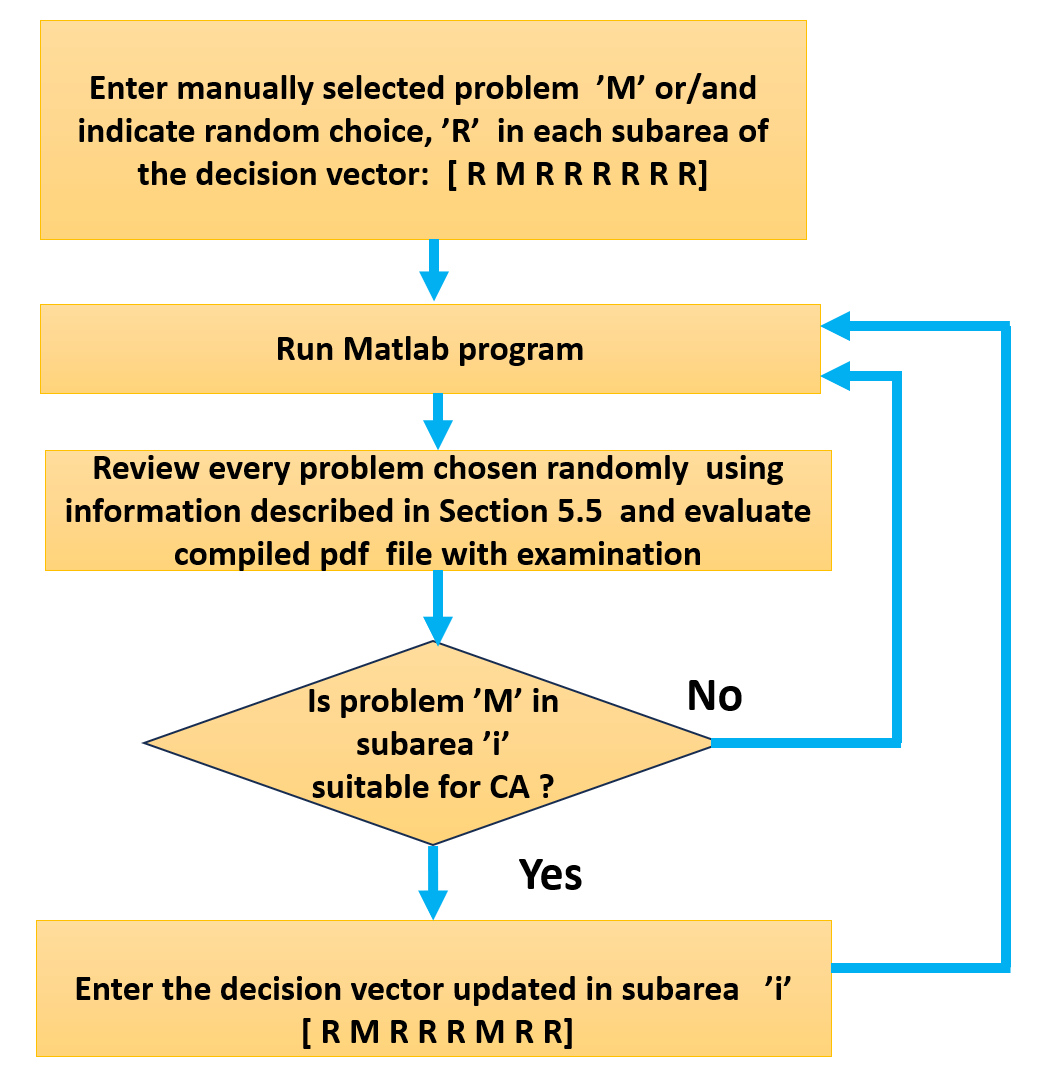}
\end{center}
\caption{\small{Flow chart for stepwise constructive alignment.
The alignment starts  with selection of an initial decision vector, where the educator
enters any manually chosen problems (M) and indicates subareas
where the problems should be chosen randomly (R). For example,
the problem which has the number $M$ is chosen manually  for the second subarea
of the decision vector, i.e. $[R~M~R~R~R~R~R~R]$, and all other problems
in the seven other subareas are chosen randomly.
\\
Educator runs the program,
reviews information provided by the system in terms of ILOs, SOLO levels, indices of difficulty etc., \cite{wace24}, evaluates compiled pdf output file and makes the decision about the inclusion of the problems in examination by updating the decision vector. The educator runs the program again with the same decision vector, if suitable problems could not be found in the trial.
The desired number of points is guaranteed by the system for each run and the system
automatically rejects the problems which have been used recently.
 After a number of reruns the educator should be able to find the best possible solution for CA with appropriate indices of difficulty.  The alignment is performed in a stepwise way where each step is associated with an update of the decision vector.}}
\label{figf}
\end{figure}

%\section{Introduction \& Main Contributions}
%\label{int}
%\noindent
%The introduction to this report starts with a concise explanation of the concept of CA (Constructive Alignment) in Section~\ref{ica} %which will be the main concept in future control education. The literature and the difficulties associated with application of this %concept are reviewed in Section~\ref{mot} and  formed the basis for further development described as contributions of this report in %Section~\ref{prs}. Organization of this report is introduced in the end of this Section.
\section{Constructive Alignment in Automatic Control Education}
\label{ica}
\noindent
Control systems are widely used in many application areas and therefore
control education is an important part of most engineering programs.
Control engineering is generally considered as a very challenging subject
from many mathematical and practical aspects.
\\
For improvement of the learning outcomes the concept of CA (Constructive alignment) developed by \cite{big} is widely applied in control courses, see for example \citeauthor{kno} \citeyearpar{kno} and  \citeauthor{sto1} \citeyearpar{sto1}. CA is based on the following principles:
\begin{itemize}
\item Educator designs the ILOs (Intended Learning
Outcomes)  using the verbs which describe the content to be
learned.
\item Educator creates the TLAs (Teaching and Learning Activities) including
the assessment methods
which address these verbs. Common verbs form the basis for constructive alignment
of teaching and assessment.
\end{itemize}
All the components such as  ILOs, TLAs and assessment tasks (as the part of learning) are based
on the same verb model in  the constructively aligned course. The results of instruction are massively improved
when curriculum and assessment methods are aligned, due to the teaching quality assessment (\citealp{big}).
The literature and the difficulties associated with application of this concept to automatic control education are reviewed in the next section,  which  formes the basis for further development described as contributions of this report in Section~\ref{prs}.
\section{Previous Work \& Motivation for Development}
\label{mot}
\noindent
Final written examinations is the most important part of summative assessment (evaluation of the  student learning at the end of the course) and the last step in the constructive alignment associated with the course.
Preparation of constructively aligned examination (which measures precisely the ILO competences) with a given number of points takes significant amount of time of the control educators, \cite{das}.
This strongly motivates development of a toolkit for automatic compilation of examination problems, which facilitates
constructive alignment. The toolkit should provide all necessary information for successful alignment of examination problems with ILOs, \cite{amr}.
\\ To the best of our knowledge, this is the first work which proposes such  a toolkit for automatic control courses. The closest work (presented by \cite{gru})  is associated with the development of automatic generation of examination problems in statistics using R.
\\ An interesting  literature review on automatic generation of examinations and the relation to  Bloom’s taxonomy can be found in \cite{nd}.
Note that the educational taxonomies, like  Bloom’s taxonomy  were designed as generic methods supposed to be the same in all subjects. However, known taxonomies are not able to fit well to all subjects, which resulted in the
development of so-called subject oriented taxonomies, for example the computer science-specific taxonomy developed by \cite{ful}.
\\
Unfortunately, the concept of CA is not directly applicable to the automatic control courses due to
possible misinterpretation of the verbs that may appear in the verb based model, which in turn results in robustness problems.
For example, the verb 'describe' can be associated with different levels of understanding.
For example, the student can 'describe' simple relationships, or complicated cascade control systems which require evaluation and analysis.
In this case different interpretations of the verb 'describe' should be recognizable by the system
and overgeneralization should be avoided (\cite{bigc}).
\\
These difficulties can be avoided via development of the SOLO (Structure of the Observed Learning Outcome) taxonomy  associated with the subject  of automatic control, and  quantification of the ILOs using the developed SOLO model (see for example \citeauthor{sto1}, \citeyear{sto1}).
Proposed solutions are outlined in the next section.
\section{Proposed Solutions \& Main Contributions}
\label{prs}
\noindent
Each learning outcome and examination problem in the database is associated with the
SOLO model, which was developed for the control course. Educator gets the opportunity  to evaluate each examination item,  (see \citealp{lio}) and to align the examination using the toolkit described in this report.
This approach is exemplified in Appendix of \cite{wace24}, where the ILOs are quantified using the SOLO model for the course
'Process control and measurement techniques'.
\\
The rapidity of the computer program (which is a measure of the software performance) is a very important requirement for this system since  multiple runs may be required for CA.
SCA (Stepwise Constructive Alignment) is a new concept introduced here and supported by the
developed toolkit. The educator reviews randomly suggested choices of problems in each step of the alignment in terms of relation to ILOs, SOLO levels, levels of difficulty, etc., and selects the problems which are relevant for CA to the next run. The method allows to align the examination step-by-step, whereby the desired number of examination points is automatically guaranteed by the program in each run.   SCA can be easily done with computer aided tools such as the decision support Matlab/\LaTeX\ toolkit developed in this report. The flow chart for SCA is presented in Figure~\ref{figf}.

\section{The Choice of Software \& Motivation for Customized Development}
\label{sp}
\noindent
The \LaTeX\ typesetting system is far superior to Word, Open Office and other text editing programs due to
\begin{itemize}
\item high typographical quality of the documents
\item automatic formatting, itemization, indexing, equation numbering and reference generation
\item flexible structures of the documents
\item  convenient typesetting of complicated mathematics, figures, tables and general technical contents
\item a large and increasing number of free add-on packages which facilitate typesetting
\item a number of customization options, which allows to create documents with user defined specifications
\item compatibility with a number of operating systems, like  Windows, Linux and MacOS
\item support from a  large number of communities like computer scientists,  mathematicians, engineers and many others
who represent both developers and users.
\end{itemize}
\LaTeX\ is widely used in the control community for preparation of collaborative research articles and other scientific documents. Also final written examinations are generally prepared in \LaTeX\ at 'Division of Systems \& Control' and
a large number of examination problems written in \LaTeX\
has been accumulated over the years. \LaTeX\ was therefore chosen in this development as the software platform  for preparation of the examination documents.
\\ One additional software platform is needed for development of the database and performing
advanced calculations required for the toolkit. Matlab is a high level programming language with flexible data structures, which are needed for development of the database.
At the same time Matlab is a well developed computational platform that includes a large number of
mathematical functions, which are needed for automatic selection of the examination problems.
Matlab comes with a number of toolboxes designed for automatic control systems and has powerful visualization tools.
\\
Many examination problems at the 'Division of Systems \& Control'  have been prepared and visualized using
Matlab. Therefore, Matlab was chosen as the second platform for development of the toolkit.
Matlab Database Toolbox  could have been used for development of the database of examination items.
However, number of limitations and the cost of this toolbox were the motivation for choosing a customized development instead.
The database was associated with nested cell arrays, dynamic structure and classification codes which address unique requirements and facilitate automatic selection for the specific course. Note that nested solutions with dynamic memory allocation are not supported by SQLite or any other standard SQL database included in the Database Toolbox.
\\
\\
Description of the toolkit is presented by \cite{wace24} :
\\ \\
Stotsky A. \& Wik T., Automatic Generation of Examinations in the Automatic Control Courses:
Decision Support  Matlab/Latex Toolkit  for Stepwise Constructive Alignment, In Proceedings of 2024 IFAC Workshop on Aerospace Control Education, July 22 - 24, 2024, Bertinoro, Italy. 
\\ IFAC PapersOnLine 58-16 (2024) pp. 47-52 
\\ \url{https://www.sciencedirect.com/science/article/pii/S2405896324012291}

\section{The Toolkit as Decision Support  System  for Stepwise Constructive Alignment}
\label{sca}
\noindent
When preparing the final examination the educator evaluates all the TLAs including formative assessments which have been
done in the course, and uses the possibilities offered by the system to align the final examination.
For example, the educator may exclude problems from subareas that have been assessed by other methods, like exercises,
laboratory works, formative tests and others. At the same time, the educator may include several problems in subareas which were not tested properly during the course. The opportunity to include manually any problem from the database is especially beneficial for CA. This option can be called CT (Choice of Teacher). The educator is able to choose several problems relevant for CA in some subareas and let the system randomly choose the rest of examination. For example, the educator may update the database and include new problems in the examination that are especially relevant for the current course.
\\
Note that each examination problem in the database is associated with the ILOs and quantified using SOLO verb model and the index of difficulty. Automatic generation of examination for a given number of points takes a few seconds in Matlab and the educator may run the program several times in order to choose the best solution.
The educator may include the problems, which randomly appeared in the previous run to the next run according to the concept of SCA (Stepwise Constructive Alignment).
In other words the educator reviews the random choice of problems in every subarea in each step, and manually selects the problems to be kept for the next run, introducing additional constraints in the next run.
The flow chart for SCA is presented in Figure~\ref{figf}.

\section{Conclusions \& Outlook}
\label{conc}
\noindent
% Constructive Alignment will be the main concept in future automatic control education.
Significant improvements of the results of learning due to application of the concept of CA
and development of modern software tools (like Matlab \& \LaTeX\ ) was the motivation for developing
a toolkit for automatic compilation of
examination problems in automatic control courses.
The decision support Matlab/\LaTeX\ toolkit, based on random number generators
for selection of examination problems, formed the basis for development
of a new strategy, Stepwise Constructive Alignment, which has the potential for further
enhancement (compared to traditional CA) of the learning.
\\ Further improvement of alignment can be associated with software development and
introduction of ILOs, SOLO and indices of difficulty in
the selection algorithm to automatically reject bad trials, choose desired level of difficulty,
facilitate SCA and provide additional support to the educator.
These improvements may result in development of a new concept of Automatic Constructive Alignment based on computer aided technologies which would allow to reduce the involvement of skilled educator/instructor (that implies significant time savings) in preparation of educational materials.
\\
Further development can also be associated with updating parameters in the examination problems  using toolboxes (like control system toolbox, filter design toolbox, signal processing toolbox, symbolic math toolbox, system identification toolbox, robust control toolbox,  and many other toolboxes) which are available in Matlab.

\end{document}